\documentclass[%
 aip,
 amsmath,amssymb,
 reprint,%
]{revtex4-1}

\usepackage{graphicx}% Include figure files
\usepackage{dcolumn}% Align table columns on decimal point
\usepackage{bm}% bold math
\usepackage[utf8]{inputenc}
\usepackage[T1]{fontenc}
\usepackage{mathptmx}
\usepackage{etoolbox}
\usepackage{tikz}
\usetikzlibrary{arrows}
\usetikzlibrary{arrows,automata}

\usepackage{natbib} %Bibliography.
\setcitestyle{numeric} %Cite as numbers or author-year.

\makeatletter
\def\@email#1#2{%
 \endgroup
 \patchcmd{\titleblock@produce}
  {\frontmatter@RRAPformat}
  {\frontmatter@RRAPformat{\produce@RRAP{*#1\href{mailto:#2}{#2}}}\frontmatter@RRAPformat}
  {}{}
}%
\makeatother
\begin{document}
\preprint{AIP/123-QED}
\title[Structure of Chaotic Attractors in FitzHugh-Rinzel model]{Large-Scale Structure of Chaotic Attractors in FitzHugh-Rinzel model}
% Force line breaks with \\
\author{Mohammadreza Razvan}
 \affiliation{Mathematics Department, Sharif University of Technology}%Lines break automatically or can be forced with \\
  \email{Razvan@sharif.ir}
  \author{Sheida Shahidi}
  \affiliation{Mathematics Department, Sharif University of Technology}

\date{\today}% It is always \today, today,

\begin{abstract}
Chaotic bursting behaviors have been observed by many authors in neural dynamics mainly in the transition between different kinds of bursting behavior. 
As a well-known three-dimensional ODEs model with various bursting solutions, the FitzHugh-Rinzel model has been considered in this paper. The structure of the strange attractor that appears in chaotic transitions of this model was investigated by introducing a stochastic approach to uncover the transition mechanism. 
To portray this idea the attractor of the dynamical system can be partitioned into some regions and a discrete evolution that is inspired by the flow between them is sketched. 
A suitable Markov chain has been associated with the strange attractor based on partition selected by the recognizable regions of the attractor. Then the entropy rate of the Markov chain and the topological entropy of dynamical systems are compared to decide if the associated Markov chain should be modified or not. Furthermore, the differences between entropies guide us to uncover some changes in the shape of the attractor including some new regions which play important roles in the chaotic behavior of our system. 
It can be also ensured with the help of Lempel-Ziv quantity that the estimated entropy rate is reliable.

\end{abstract}

\maketitle

\begin{quotation}
The stochastic approach is widely used in the study of chaotic dynamics. There are  several stochastic tools in dynamical systems such as Shannon entropy,  Lyapunov exponents, Markov transition matrix, Kneading sequence, Kaplan-York dimension, etc. In classical Markov transition matrix approach the attractor is carefully divided into some regions and the dynamical systems is reduced into shift dynamics between different regions. Then one tries to find an equivalence between the original dynamics and shift dynamics on the associated Markov chain. However, it is not always easy to find such a rigorous equivalence but some attractors have specific zoning that facilitates this approach. For instance the Lorenz attractor could be divided into left and right regions. Other examples are in one dimensional dynamics in which kneading sequence approach is applicable. In this paper we are confronted with the chaotic attractors of FitzHugh-Rinzel that is portrayed as a Leviathan with the algorithm of beheading and creating new heads and wings.
\end{quotation}

\section{Introduction}

%\textit{Importance of Choatic Bursting}\\
As with any dynamical system, various behaviors can be seen in neural models such as continuous spiking, bursting, and mixed-mode oscillations(MMOs)\cite{DELNEGRO1998174,Rinzel1987}. Chaotic mode is a prevalent performer to initiate, transform, and terminate these behaviors in neural models\cite{Barrio2012Kneadings,Medvedev2006Transition,Shw2011,10.1063/1.166488,Terman1991}. It also appeared as an intermediate state for transition between different types of these oscillations\cite{CHAY1985357}.
In this paper chaotic bursting in the FitzHugh-Rinzel(FHR) model as a neural model has been investigated. Likewise other research results chaotic attractors appear in the transition between bursting behaviors. The attractor creation, change, and disappearance are presumed to be vital factors in the evolution of neural dynamics.

%\textit{Discretization}\\
In order to investigate chaotic attractor, it is carefully partitioned into some seemingly independent regions to define states of the Markov chain.
 A directed graph associated with the dynamic is created by the most probable state that is dictated by the flow. In other words the random walk on the graph is induced by chaotic dynamics. The structure of this graph will be used to estimate the related stochastic quantities of the system to investigate the complexity of the dynamics.
 This idea leads to a rough description of a chaotic system as a finite-dimensional stochastic process. 

%\textit{Stochastic Approach}\\

The stochastic approach is widely used in dealing with chaotic dynamics and its features are helpful to study the chaos in dynamical systems. Several papers used stochastic tools such as Shannon entropy,  Lyapunov exponents, Markov transition matrix, Kneading sequence, Kaplan-York dimension, and Lempel-Ziv complexity.
These measures apply to a wide range of problems in science\cite{Barrio2012Kneadings,Paun2020,doi:10.1142/S021812742230004X}. Recently, in \onlinecite{Scully2021Measuring}, Scully et al. applied these computational indicators to Lorenz and Rossler attractor to measure the chaos.
Here the entropy rate of the Markov chain is computed and it is compared with the topological entropy of dynamical systems to estimate the complexity of this random walk. Topological entropy is a usual way to measure the complexity of dynamical system by viewing its distinguished orbits growth rate\cite{Brin2002Introduction}. The ability of the Lemple-Ziv sequence is also applied to ensure that the above walking is random\cite{Lempel1976}. Indeed, the Lemple-Ziv complexity can be used as a correction term for the entropy rate.

  %\textit{Information Theory}
  This method is indebted to I.Tsuda by his idea of chaotic itinerancy\cite{Tsuda:2013,doi:10.1063/1.1607783,Tsuda2001,Tsuda1987}. His great work for modeling episodic memory introduced regions resulted from attractor ruins that the flow itinerant between them. The dynamics spends a long time in the vicinity of each region in the phase space. In this paper an attractor ruins so many times instead of ruining so many attractors suddenly. The impact of each region of chaotic behavior could be individually examined in this way. 
  
There is an excellent systematic approach that is defined by the Frobenius-Perron operator to perform dynamical coarse-graining for low dimensional chaotic systems. This is an information theory approach to chaos by considering the propagation of a sequence of the letters chosen from the appropriate alphabet with the advantage of using the stochastic prospect\cite{Nicolis2015Chaos,Basios2011Symbolic,MacKernan2009LocalGlobal}. 

One can consider finitely many states on the attractor and construct the Markov Model associated to the dynamics of the attractor. This can also be done for a time series by applying Taken's embedding approach\cite{hou2015dynamical}. This Markov model is known as Ordinal networks introduced by Small in 2013\cite{small2013complex}. A multi-step  HMM is apparently more efficient at extraction of deterministic and stochastic parts of the dynamics from a given time series. 
%In some sense A suitable VAEs could be optimal to achieve this goal.
The main goal of time series analysis is prediction but we are more concerned with the extraction of geometric features of chaotic attractors which is quite meaningful in Tsuda's view point to memory. 

%The recurrence network introduced by Donner et. al in 2010 and used recurrence in phase space to construct adjacency matrix of an associated complex network. 
%The ordinal network approaches begin with the work of small et. al in 2017 and mapped time series to the set of ordinal symbols and build a network with the aid of symbolic dynamics.
%Theses approaches suffered from  lack of distinguishing between chaotic dynamics and random dynamics. 
%These ideas with the help of Taken's theorem have been used to developed the concepts of recurrence networks and ordinal networks. These approaches used the  complex network  to obtain some dynamical properties of times series. 

There are several ways that a continuous dynamical system can be coarse-grained. Here the anatomy of Leviathan has been used to discretize the fine-grained dynamical system which is a more intuitive way by getting involved in the nature of chaotic attractors in the FH-R model.
 
  %\textit{Leviathan}
  It has been tried to apply this method to study the geometry of attractors via the change in the parameter. The geometry of the attractors is evolved by appearing some new parts so called heads and hands.We called such an attractor the "Leviathan" because of its several heads and hands. It appears that more heads and hands facilitate the attractor to be more chaotic.  The anatomy of these Leviathans was used to estimate the related stochastic quantities of the system for investigating the complexity of dynamics.\\
  In this framework, we associate a directed graph to the Leviathan with a coarse viewing of the dynamics of a chaotic attractor.
  Whenever a new head or wing appears a new vertex is added to the graph. 
  This approach is coarse-grained since it only considers appearing or disappearing and merging some heads and wings or dividing some regions. 
  The information concealed in this coarse-grained approach helps us to understand the underlying dynamics.
  
 %\textit{Entropy}
 In this paper the concept of entropy is applied to describe large-scale changes in chaotic attractors. Entropy can help our large-scale viewpoint of the chaotic attractor because of its nature in describing the long-term behavior of the stochastic process. The information that is carried in this way is the concept that is common in the dynamical systems and random processes by the work of Shannon\cite{Shannon1948MathematicalTheoryII}. This information and its changes are what is employed to investigate the structure of the chaotic attractor.
 In order to estimate the complexity of this random walk the entropy rate is computed. Since this random walk has been extracted from a chaotic dynamical system it is reasonable to compare this entropy rate with the topological entropy of dynamical systems. Indeed the Lemple-Ziv complexity can be used as a correction term to the entropy rate.
 By the nature of this coarse view, indubitably some changes in the attractor are missed and it tries to find them by attending to the difference between these two types of entropy as an indicator for information loss. As it will be described later in details, by this approach some new hands and heads and other large-scale changes in the structure of the attractor could be detected.

This paper is organized as follows: first in section II, the FH-R model was introduced in details by the historical way of neural models. Section III is devoted to describing the method.

\section{Neural Model}
The FitzHugh-Rinzel model, a reduction of the Hodgkin-Huxley(H-H) model introduced by Rinzel in 1987, is used as a neural model. 
Hodgkin and Huxley have modeled the propagation of an action potential in the squid's axon in 1952 by a four-dimensional ordinary differential equation \cite{Hodgkin1952Quantitative} and the neural models entered the arena.
Subsequently, some variations have appeared for the aim of better understanding the features, capacity and reliability of H-H model that one of them is FH-R model. 

FH-R model is a well-known elliptic burster derived by adding a linear equation to FitzHugh-Nagumo(FH-N) model, as a slow variable\cite{doi:10.1142/S0218127400000840}. FH-N is a two dimensional reduction of H-H that was attained in the effort of solving H-H in 1961 by using the van der Pol equation\cite{Rinzel1987}. 
This model was called the Bonhoeffer-van der Pol model by Rinzel. In 1962, a circuit was associated by Nagumo to this system of equations so it was called the Fitzhugh-Nagumo model. 

The mathematical model of FH-R which is known as a reduction of the Hodgkin-Huxley(H-H) model for neurons is given by the following three-dimensional systems of equations:
\begin{center}
	$v' = v - v^3 - w + y + I$\\
	$w' = \phi (v + a - bw)$\\
	$y' = \epsilon(-v + c - dy)$
\end{center}
In the original paper of Rinzel, parameters were considered as $I=0.3125$, $\epsilon = 0.0001$, $\phi = 0.08$, $a=0.7$ and $c=-0.775$.
%Yadav et. al. have discussed the bifurcation analysis of FH-R and the evidence of chaotic bursting by Lyapunov exponent analysis.[Ref. Yadav] (can't fine reference!)
In \onlinecite{Wang2021} Wang et al. discussed the FH-R with the impact of magnetic field. They also considered a network of FH-R model and studied the existence and propagation of wave solution affected by external stimuli.
Because of the fast-slow nature of the FH-R system \cite{Rinzel1987}, applying the singular perturbation methods may be useful to understand the behavior of the system. For example, Xie et. al. use this method to find some so-called double MMOs and double canards\cite{XIE2018322}.
By varying some parameters FH-R has also been considered as a three time scale system and the MMO of that system has been studied\cite{doi:10.1137/100791233}.

FH-R is a reach and well-known model that has continuous spiking mode, bursting, chaotic behavior, and MMO. The transition between these regimes has been considered in some studies\cite{doi:10.1063/1.5042078}. In Shilnikov 2011, authors discussed MMO phase in the transition between bursting and quiescence by using the Voltage Interval Mapping technique\cite{Shw2011}. They also found a torus attractor in FH-R and applied the kneading sequence approach to approximate entropy in the chaotic regions. It is well known that the role of chaos is substantial in these transitions and  it is reasonable to focus on theses chaotic attractors. 

 In 1998 Del Negro et al. introduced an equivalent system of equations\cite{DELNEGRO1998174}:
\begin{equation}\nonumber
\begin{matrix}
v' = aw - 4v^3 + 4v -z \\
w' = -(1 + 4v + w) \\
z'= \alpha(bv - (cz- z0)/d)
\end{matrix}
\end{equation}
with  $a=1$ ,$\alpha=0.003$ , $b=1.25$ , $c = 1$ , $z_0=1.33$ and $d=4$. In contrast to the original FH-R model, it has two fast and one slow variables.
In this study, "$a$" has been taken as a parameter and the other one fixed as $\alpha =0.006, b=6, c = 1.605, z_0=1.1$ and, $d=3.7$. With these choices various bursting and MMOs would be seen in this model.\\
One of the interesting features of this model is the existence of complicated  bursting attractors in the range of parameter $a=0.7139$ to $a=0.718$.
Bursting is a prevalent oscillation in neural models\cite{Su4173,Womack10603,Pace2007Inspiratory,Izhikevich:2006}.
It is well known as a mechanism for convey information between different parts of nervous system that reflect neural coding\cite{SHERMAN2001122,Izhikevich2003Bursts}. It is a periodic behavior that can be recognized by grouped bursts of spikes following by a period of quiescence.   
Almost all mathematical models for busting neurons use slow-fast dynamics to described it. FHR is a famous model which exhibits various bursting behaviors by varying the parameters in this reduced neural model.
%In the range of parameters that the system was studied, three fixed points exist for the slow subsystem that the system is surely affected by the bifurcations underlying them.  We are aware that the change in the attractor especially in the non-chaotic region can be understood by studying the bifurcations of the slow system and the slow-fast dissection is important in the study of bursting solutions. 
%In this study the importance of the slow-fast dissection was dismissed and 
The attractor has been viewed as a complex oscillator, and the geometric structure of it is crucial in our study. For this aim, helicopter view of the attractor has been used, main features of geometric structure extracted and major changes of them detected.

We start with a simple periodic orbit at the parameter a=0.7138  that can be seen in Fig. 1 and progress by changing the parameter 
\begin{figure}[h!]
	\begin{center}
		\includegraphics[scale=0.185]{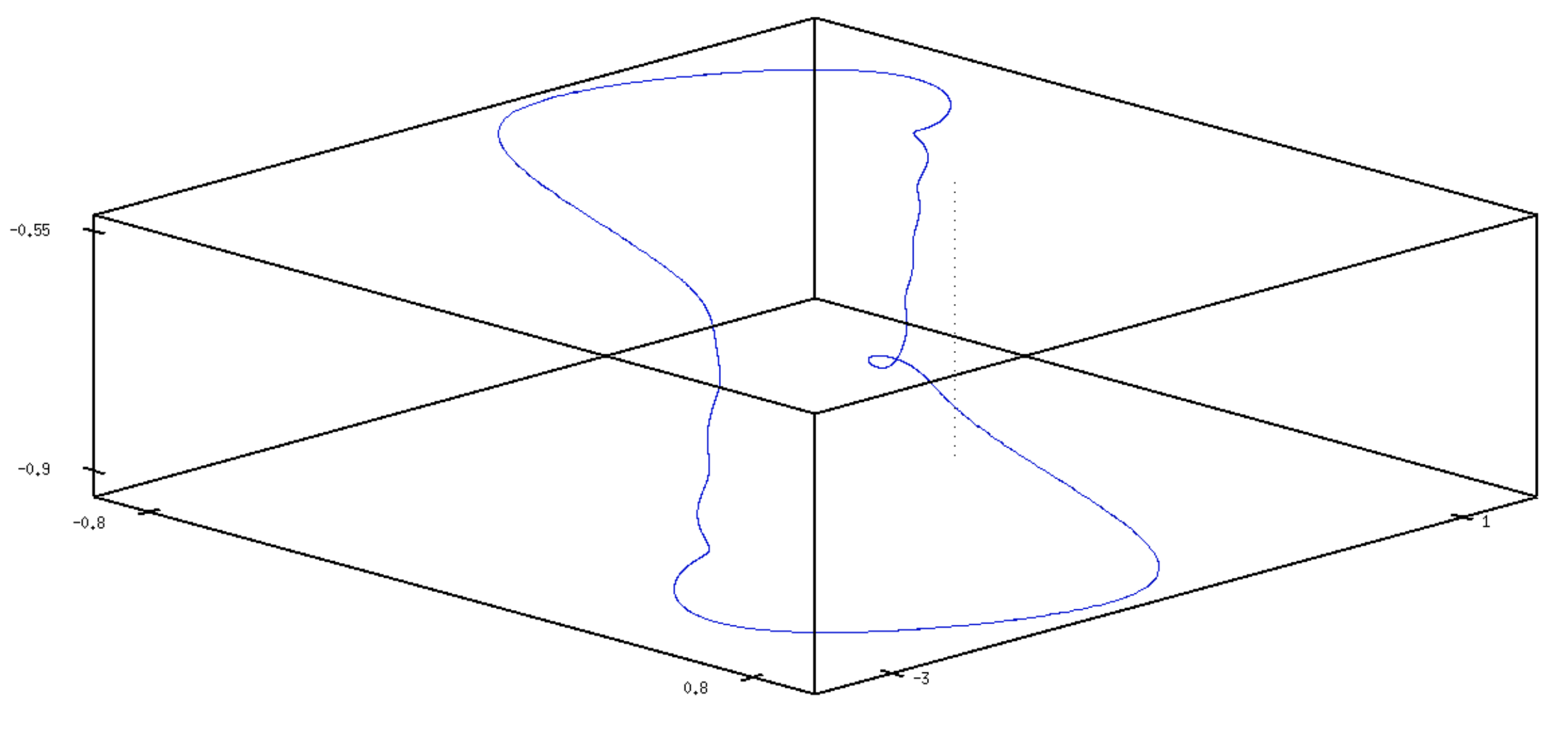}
	\end{center}
	\caption{Simple Periodic Orbit}
\end{figure}
At the end of the study it become to a two layer periodic orbit at parameter a=0.7178 (Fig. 2).
\begin{figure}[h!]
	\begin{center}
		\includegraphics[scale=0.185]{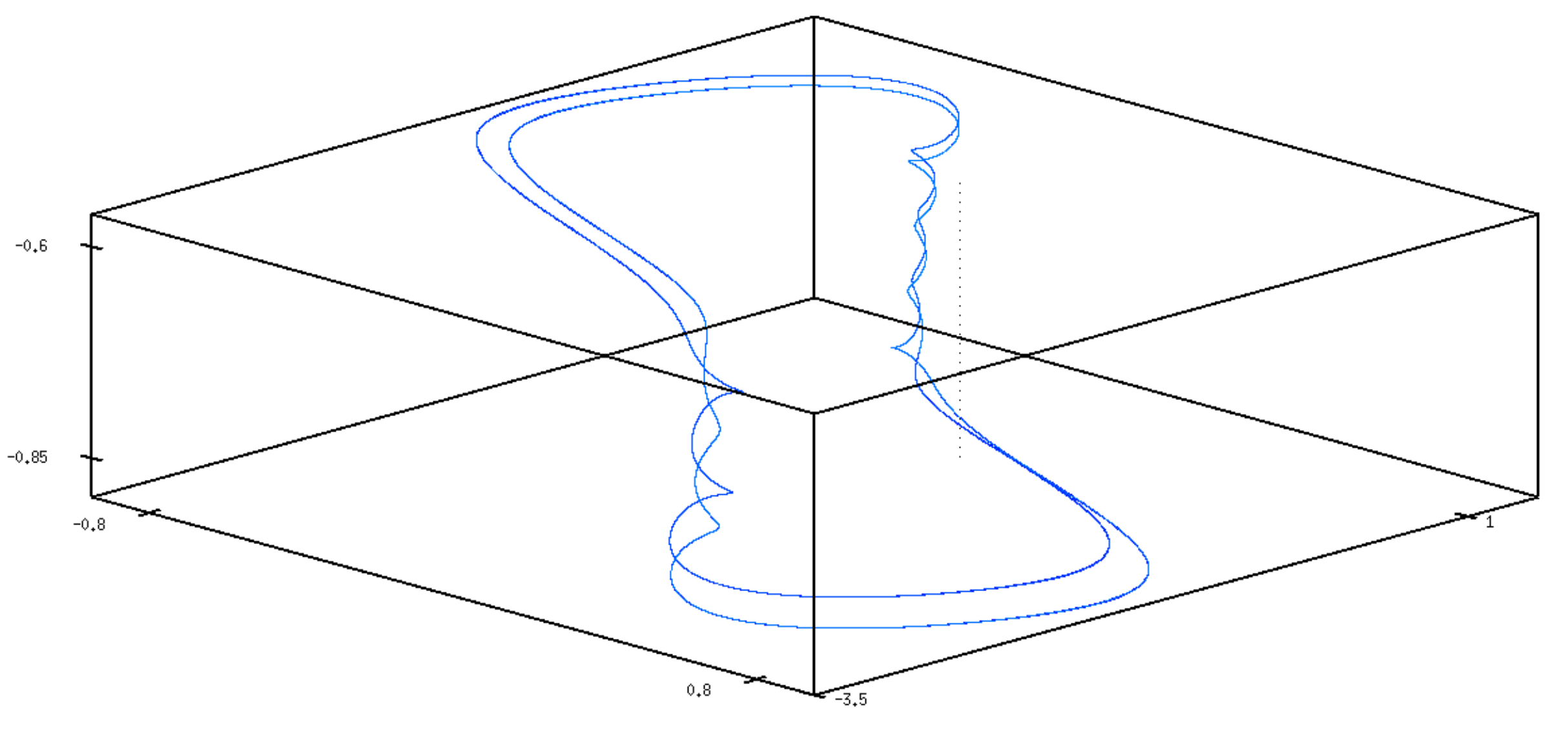}
	\end{center}
		\caption{Two Layer Periodic Orbit}
\end{figure}
There are two periodic orbits in this range of parameters. By choosing a suitable Poincare section, chaotic regimes can be verified in the system.(Fig. 3)
\begin{figure}[h!]
	\begin{center}
		\includegraphics[scale=0.185]{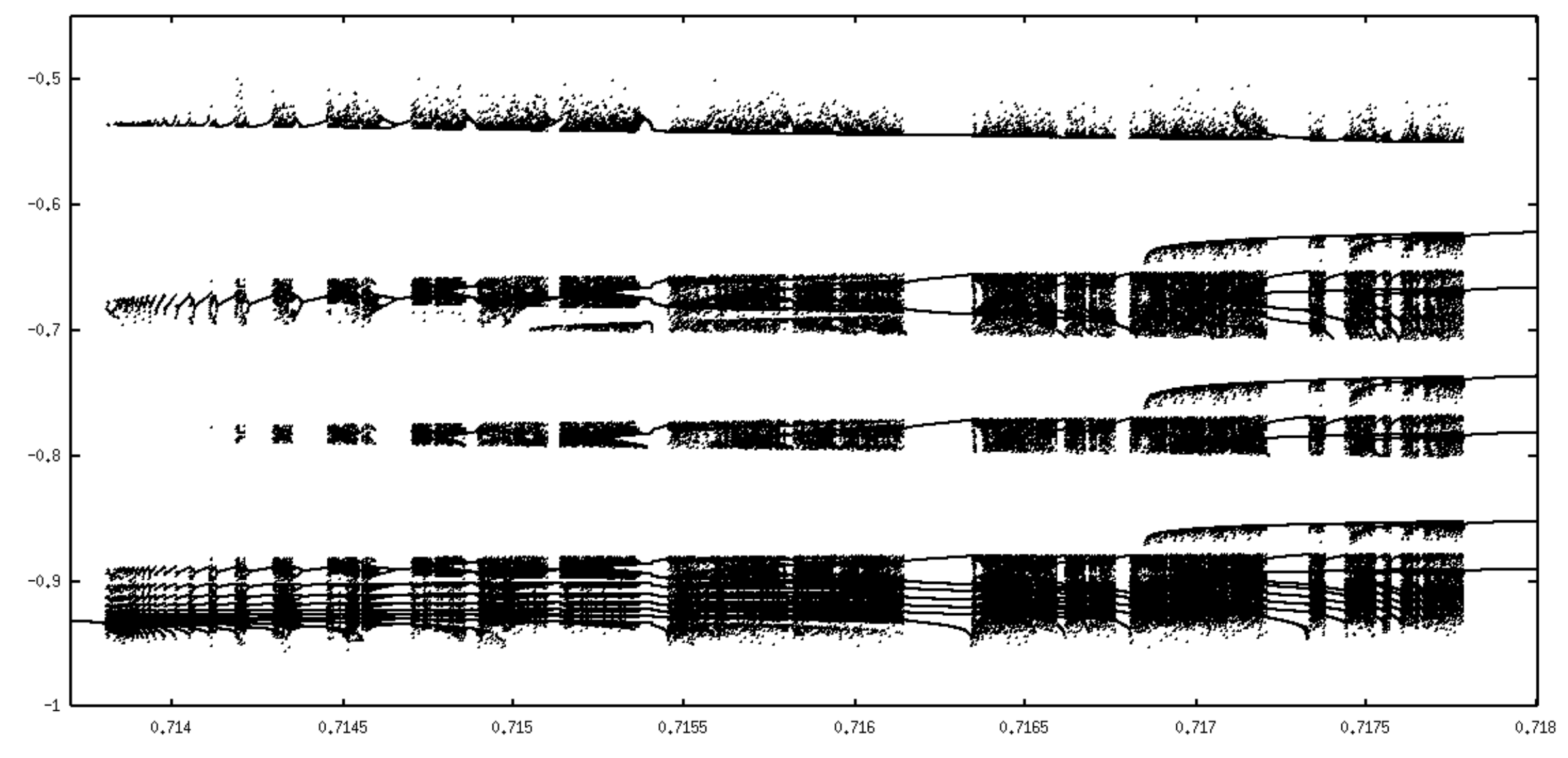}
	\end{center}
		\caption{Poincare map in range of parameter 0.7136 < a < 0.718 }
\end{figure}
By varying the parameter various types of Leviathan appear witch will be investigated in the next section.

\section{Method}
The current research utilized stochastic techniques to extract dynamical features using coarse-graining approach.
Some important regions on the attractor are selected as "Markov partitions", and a Markov chain is obtained by following the dynamics. A random walk can be taken on the graph that is attributed to the Markov chain as the dynamics evolves.
We have a chaotic dynamical system on one hand and a Markov chain extracted from it on the other hand. 
When the Markov chain is obtained by using the prospect of information theory, entropy rate as one of its features is compared to entropy of the primary dynamical system to verify the validity of the Markov chain and randomness of the flow defined walking.
To use the Perron-Frobenius theorem, the derived Markov chain is assumed to be aperiodic and irreducible, and the entropy rate of the associated random walk can be computed.

The attribution  of Markov chain is a prevalent way to discretized a continuous dynamical system.
When a chaotic system is studied, the underlying continuous dynamical system is converted to a discrete one with the help of symbolic dynamics.
The Poincaré mapping theory can be utilized by considering a suitable section. The section is divided into finite regions and the dynamics is transformed into the symbolic dynamics with the number of codes equal to the number of the regions. 

In a more general situation, phase space can be divided into a finite set of cells, and the statistics induced by deterministic dynamics becomes the problem under consideration.
Thus, a Markov chain is constructed from the chaotic dynamics whose rules obey the laws of the non-linear dynamical system\cite{Nicolis_1986}. 

A graph can be induced by the associated Markov chain and a random walk on it is inspired by the coarse-grained chaotic dynamics. For this aim, vertices are associated with the most recognizable region in the attractor and edges of the graph are flow defined naturally. 

The attribution of graphs to dynamical systems has been considered in some research with different mechanisms to find features of the associated random walk and better understand the dynamics\cite{Nagaraj2017,Basios2011Symbolic,Perkins2014}.

For example, as an obvious choice, a directed graph with two vertices can be associated with the Lorenz attractor\cite{Scully2021Measuring}. Each wing in the attractor is associated with a vertex and the weight of edges is induced by the dynamics. The obtained graph is a complete directed graph: 
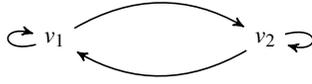
\begin{figure}[h!]
\begin{center}
	\begin{tikzpicture}[->,>=stealth',shorten >=1pt,auto,node distance=2.8cm,
	semithick]
	\tikzstyle{every state}=[fill=blue,draw=none,text=white]
	
	\node[]         (B)  {$v_1$};
	\node[]         (C) [right of=B] {$v_2$};
	
	\path   (B)  edge [bend left] node{}(C)
	(C) edge [bend left] node{} (B)
	(B) edge [loop left] node {} (B)
	(C) edge [loop right] node {} (C);
	\end{tikzpicture}
	\caption{Markov Chain of Lorenz Attractor}
\end{center}
\end{figure} 
The graph weight of this random walk will change by varying the parameters. The attractor's symmetry can be broken because of the existence of two bistable attractors in some range of parameters.
Also, a graph with two vertices can be associated with the Rossler attractor by considering a threshold on one of the variables\cite{Scully2021Measuring}.
%Another way to extract a Markov chain from a dynamical system is by Mac Kernan and Basios [2011]. Each vertex is proposed as a mean of one variable and the probabilities in the Markov chain are obtained from the propagation of time series.
%(Ref I.Tsuda). His Idea results in an attractor with some regions in the phase space that the trajectories itinerants between them.  In each of these regions, the Poincaré section is selected to discretize the system and the Markov chain will be obtained in a natural way.
%The Markov transition matrix that is obtained for the usual parameter $r=28 , s=10$ and $b=\dfrac{8}{3}$ of the Lorenz attractor can be considered as follow
%$$
%\begin{bmatrix}
%0.57 &0.43 \\
%	0.37& 0.63\\
%\end{bmatrix} 
%$$
%. For example, one different transition matrix can be obtained by changing the $r$ parameter from $28$ to $29$. 
%$$
%\begin{bmatrix}
%	0.59& 0.41 \\
%	0.49 & 0.51\\
%		\end{bmatrix}
%		$$
%But this delicate change can not be recognized easily from the shape of the attractor
%		\begin{figure}[h!]
%	\begin{center}
%		\includegraphics[scale=0.2]{Lorena29.png}
%	\end{center}
%\end{figure}   
%To the Lorenz attractor with usual parameter $r=28 , s=10$ and $b=\dfrac{8}{3}$%  the shape of attractor is like this:
%%   \begin{figure}[h!]
%%	\begin{center}
%%		\includegraphics[scale=0.7]{lorenz11.pdf}
%%		\caption{Entropy Rate}
%%	\end{center}
%%\end{figure} 
\subsection{Graph Attribution}
%(I think that we "Attribute" a Graph to a Leviathan because there is a causal relationship between regions of the attractor and vertices of the Graph!)\\
The FitzHugh-Rinzel model with coefficient under study, has a Leviathan with a recognizable head as shown in Fig. 1. The head can be seen in the middle right of the attractor. 
  \begin{figure}[h!]
	\begin{center}
		\includegraphics[scale=0.2]{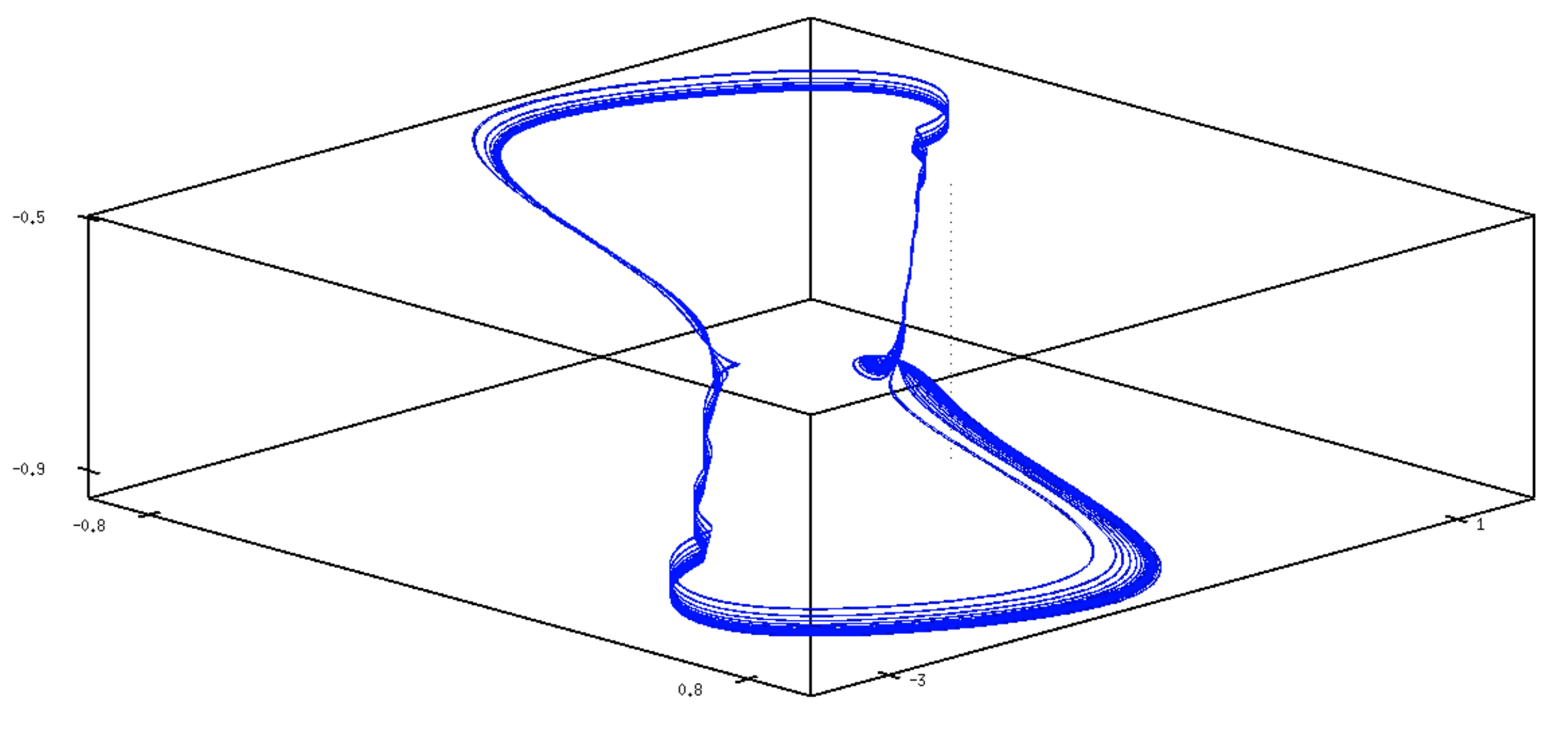}
		\caption{primitive Leviathan at a = 71385}
	\end{center}
\end{figure}
Size, structure, and number of heads evolves by parameter changes.
The creation of new heads or wings of the Leviathan reflects more information of the complex system. The appearance or vanishing of chaos can be persuaded by these new heads and wings.

The entropy is a well known tool for quantification of system complexity and uncertainty. Since the complexity of Leviathan is assumed to be summarized in the number of heads and wings the entropy changes are utilized to detect new ones of them, when they are unforeseen. 

Similar to the procedure introduced in the previous subsection for the Lorenz attractor one can attribute a graph to any Leviathan.
The vertices of graph consist of heads, wings and other important parts of attractor. The directed edges of the graph are flow defined. Some vertices with only one input and output will be neglected so far as the structure of the Leviathan is preserved. Initially a simple graph with three vertices can be attributed to the Leviathan in Fig. 5.
 \begin{figure}[h!]
 \begin{center}
\begin{tikzpicture}[->,>=stealth',shorten >=1pt,auto,node distance=2.4cm,
                    semithick]
  \tikzstyle{every state}=[fill=blue,draw=none,text=white]

   \node[]         (B)  {$v_1$};
  \node[]         (C) [below right of=B] {$v_2$};
  \node[]         (D) [below left of=C] {$v_3$};
  
  \path   (B) edge [bend left] node {} (D)
                 edge [bend left] node{}(C)
           (C) edge [bend left] node{} (D)
           (D) edge [bend left] node {} (B);
              \end{tikzpicture}
              \caption{The Graph Of Primitive Leviathan}
\end{center}
\end{figure}
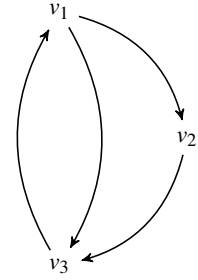
 One vertex for the top region of the Leviathan, one for the bottom, and one for the head that appears in the middle right. 
It is a primarily graph associated to the Leviathan at the first glance sufficient to begin the procedure. With the help of entropy, the graph can be refined so that it better explains the dynamics of the Leviathan which is discussed in the next section.
\subsection{Entropy Rate vs. Topological Entropy }
We introduced the walking on the attributed graph by running the dynamics.
We start from a point in a region of the attractor which corresponds to a vertex of the graph and then follow the dynamics. When another region is visited, a path from one vertex to another on the graph is chosen.
In the coarse-graining approach, when the dynamics is chaotic there is not a unique vertex could be visited in each step i.e. there is a randomness in the walking. One can compute the entropy rate of this random walk while the parameter varies.
 \begin{figure}[h!]
	\begin{center}
		\includegraphics[scale=0.5]{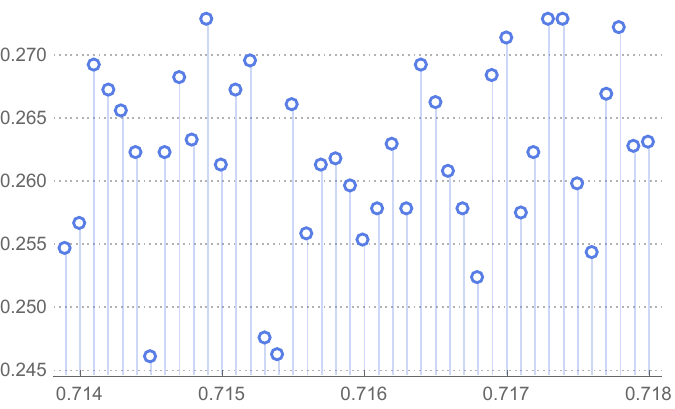}
		\caption{Entropy Rate in Terms of Parameter a}
	\end{center}
\end{figure} 
On the other hand, one can also compute the topological entropy of the Leviathan as an attractor of the dynamical systems.(Fig. 8)
\begin{figure}[h!]
	\begin{center}
		\includegraphics[scale=0.5]{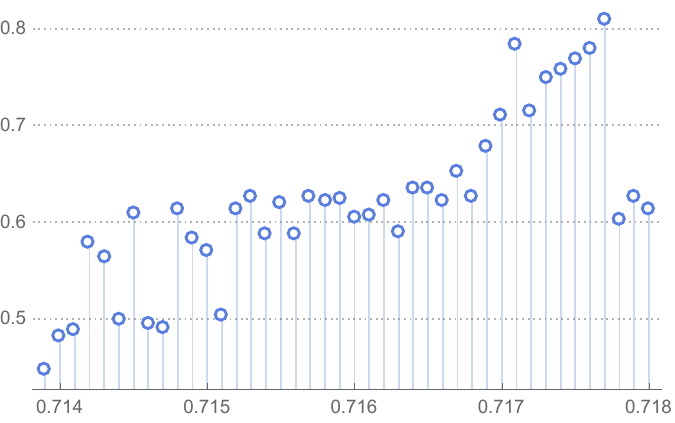}
		\caption{Topological Entropy in Terms of Parameter a}
	\end{center}
\end{figure}
There are two values of entropy. One of them is naturally the topological entropy of main dynamical system. Another way to extract the stochastic soul of the chaotic dynamics is the entropy rate of the attributed Markov chain.
There is a noticeable difference between these two values of entropy, which is not surprising because the attributed graph is very simple.
The smaller disagreements give hope that the obtained random walk has extracted more complexity of th   e dynamics. The result of comparison between these two entropy is shown in Fig. 9:
 \begin{figure}[h!]
	\begin{center}
		\includegraphics[scale=0.6]{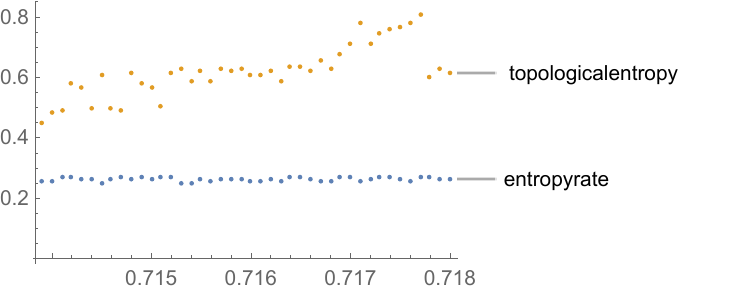}
		\caption{Comparision of Entropy Rate and Topological Entropy}
	\end{center}
\end{figure} 
\subsection{Changes in Leviathan}
 A significant difference between entropy rate and topological entropy means that the attributed graph to the dynamical system can not extract the essence of its dynamics. Also, because the entropy rate is obtained from the random walk on the graph, the difference demonstrates that the random walk can not draw out the whole chaotic dynamics of Leviathan. This is possibly because some important regions in the Leviathan are neglected in the graph extraction procedure that have a key role in dynamics. 
 The major difference between these two entropies is at the parameter value $a=0.7175$. It can be seen that a new head for the Leviathan has appeared in a small neighborhood of this value of parameter as shown in Fig. 10.
 \begin{figure}[h!]
	\begin{center}
		\includegraphics[scale=0.2]{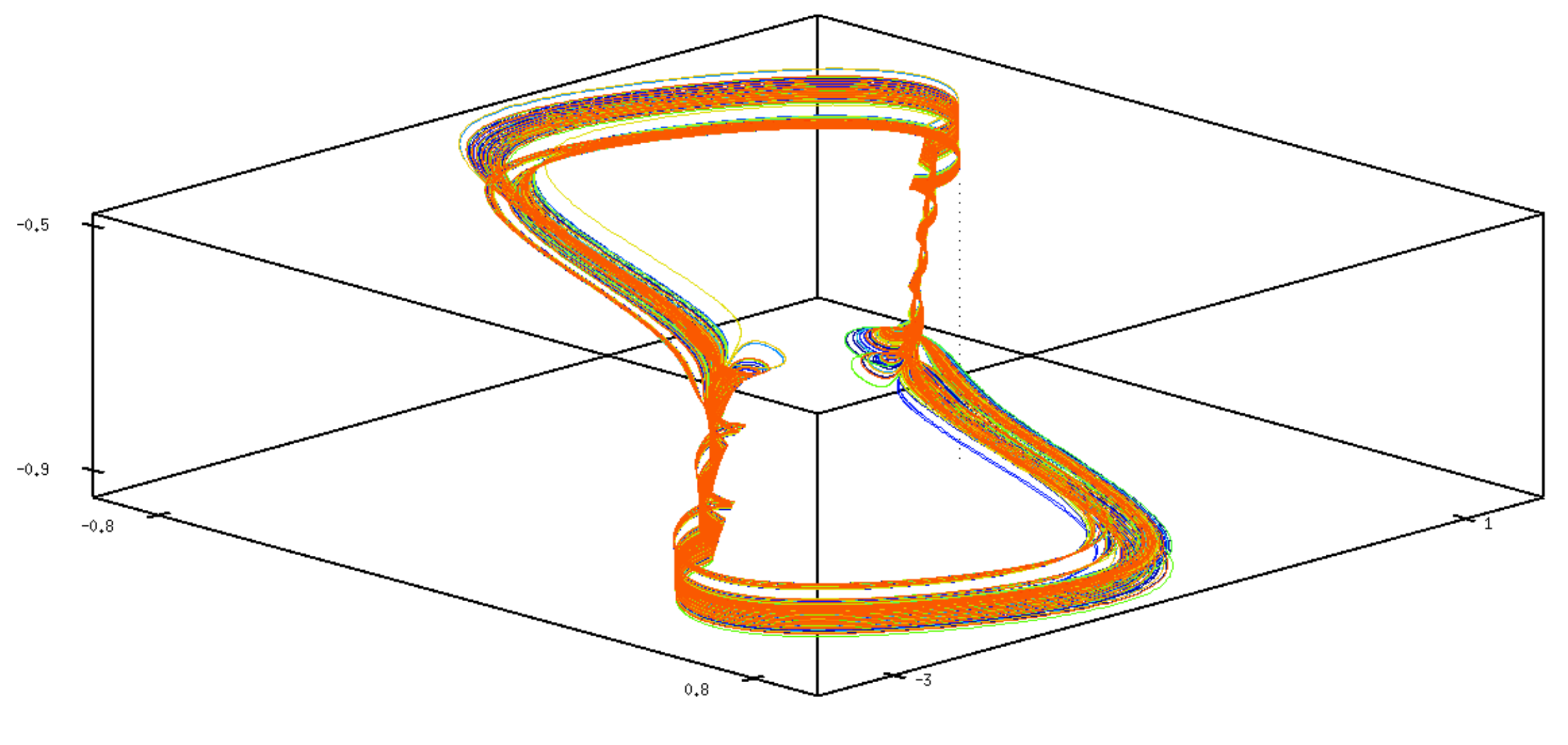}
		\caption{An Advanced Leviathan}
	\end{center}
\end{figure} 
The change in the dynamics caused by this new head also suggests the necessity of adding a vertex to the attributed graph. 
%As a result of this new head the precision of the Markov chain extraction is increased.
%It is expected that by the better graph associated to the Leviathan, we can see the birth of some new heads in the %Leviathan or the beheading of it in the specific parameters. 
 we can proceed by the new Markov chain associated to the Leviathan and it is expected to  see the birth of some new heads in the Leviathan or the beheading of it. 
 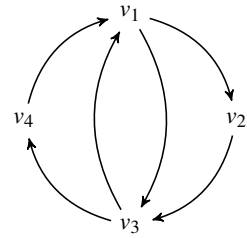
\begin{figure}[h]
\begin{center}
\begin{tikzpicture}[->,>=stealth',shorten >=1pt,auto,node distance=2cm,
                    semithick]
  \tikzstyle{every state}=[fill=blue,draw=none,text=white]

  \node[] (A)                    {$v_4$};
  \node[]         (B) [above right of=A] {$v_1$};
  \node[]         (D) [below right of=A] {$v_3$};
  \node[]         (C) [below right of=B] {$v_2$};

  \path (A) edge [bend left] node {} (B)
           (B) edge [bend left] node {} (D)
                 edge [bend left] node{}(C)
           (C) edge [bend left] node{} (D)
           (D) edge [bend left] node {} (A)
            edge [bend left]  node {} (B);
      \end{tikzpicture}
         \caption{The Graph of Advanced Leviathan}
\end{center}
\end{figure}
A new vertex for the associated graph is suggested by the new head of Leviathan and as a result a new random walk will be suggested by the modified graph that was shown in Fig. 11. 
 The modified graph better describes the chaotic dynamics. It can be verified by computing the entropy rate of the random walk on the new graph and comparing topological entropies.
\begin{figure}[h!]
	\begin{center}
		\includegraphics[scale=0.5]{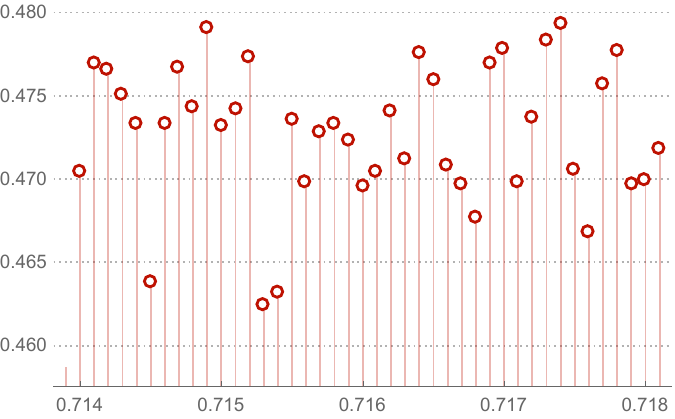}
		\caption{Modified Entropy Rate in Terms of Parameter a}
	\end{center}
\end{figure}
As it is indicated by Fig. 6 entropy rate increases for the new graph and the difference between the magnitude of entropy rate and topological entropy decreases, which shows the efficiency of graph modification. This method can be applied to any noticeable differences between topological entropy and entropy rate and figure out the hidden change in the chaotic Leviathan by finding a new head and also beheading. This procedure can be repeated so far to find the delicate changes that just a glance at the attractor can not observe.
\subsection{Lemple-Ziv}
The entropy rate of associated random walk suffers from  some overestimation when the walking is not completely random.
Lemple-Ziv Complexity is a standard tool that is applied to overcome this problem. To simplify the computations, the graph is reduced to a simpler one that is shown in Fig. 13. The vertex $v_3$ is eliminated, and a loop is added instead of it(compare Fig. 6 and Fig. 11).
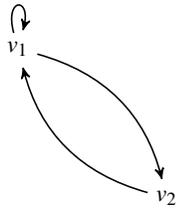
\begin{figure}[h!]
\begin{center}
\begin{tikzpicture}[->,>=stealth',shorten >=1pt,auto,node distance=2.8cm,
                    semithick]
  \tikzstyle{every state}=[fill=blue,draw=none,text=white]

   \node[]         (B)  {$v_1$};
  \node[]         (C) [below right of=B] {$v_2$};
  
  \path   (B)  edge [bend left] node{}(C)
           (C) edge [bend left] node{} (B)
           (B) edge [loop above] node {} (B);
              \end{tikzpicture}
\end{center}
		\caption{Reduced Graph of the Primitive Leviathan}
\end{figure}
A binary sequence is obtained by walking on the new graph that helps to apply Lempel-Ziv to ensure that the flow-defined walking is random. Lempel-Ziv is computed for a wider range of parameters as it is shown in Fig. 7 where in the beginning and at the end it is known that there is no chaotic behavior to emphasize the ability of this tool.
\begin{figure}[h!]
	\begin{center}
		\includegraphics[scale=0.5]{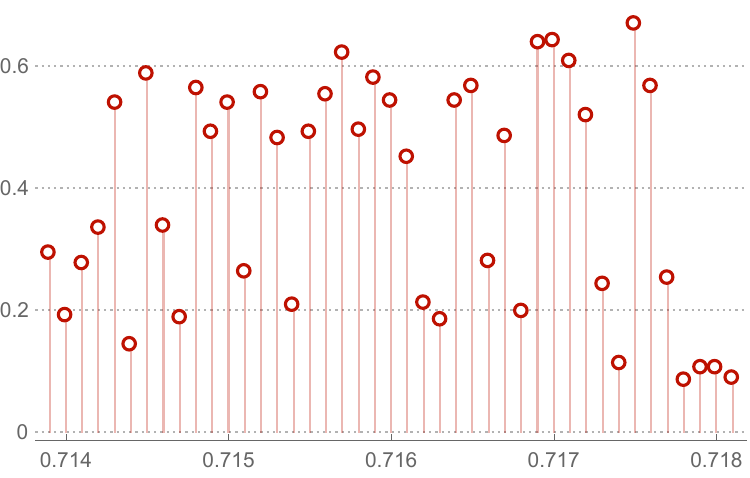}
		\caption{Lemple-Ziv in Terms of Parameter a}
	\end{center}
\end{figure}
Small values of Lemple-Ziv can be a guide to the periodic orbits with large periods.
\section{Discussion}
When dealing with complex systems, considering the system's elements is important. The state of the whole system can be influenced by the interactions of its parts. Consideration of the nature, geometry, and dynamics of this parts is necessary for understanding the individual roles that must be considered in the interaction between elements in the complex system.

A brain is a complex system in our cognition, and the complexity of its element in any layer is unavoidable. Oversimplification of elements in complicated systems can be outwitted and may conceal important mechanisms.  Neurons are elements of this system and their complexity is significant in the investigation of neural systems. The intelligence concealed in neurons of any part of the central and peripheral nervous system can change the overall neural dynamics.
The complexity of neural system is not only from the interaction of neurons but also influenced by the individual parts of this complicated system.

In the history of neuroscience, the effort of modeling neural dynamics leads to simple low-dimensional models that can be described by simple ordinary differential equations. In this paper the ability of the three-dimensional FHR model, that can describe satisfactorily the characteristics of neurons by a low-dimensional dynamical system, as the part of complicated system has been used. The structure of attractors in the FHR model will be discretized and the concentration of the  flow in different parts of each attractor has been considered.
This idea is indebted by a couple of previous significant research works, the itinerant approach of investigating dynamical systems introduced by I.Tsuda and the coarse-grained deal of dynamical systems of J.Nicolis. The advantage of these two approaches together with the grant of stochastic formalism in the realm of the coarse graining view enables us to study chaotic dynamics of neural model.

 The considered chaotic attractor of FHR model is linked to the Leviathan in our terminology because it seems to be the creature with various heads and wings that is associated to the different regions of attractor. In this playground a Markov chain is extracted from the dynamics, then capacities of the stochastic tool can be employed. It is tried to use information theoretic tools like entropy as a reciprocal concept between dynamical systems and stochastic systems, to understand the complexity of FHR model of neuron. The defect of associated Markov partition can be disclosed by detecting differences between the topological entropy of dynamical systems and entropy of Markov chain. Modification of Markov chain with the appearance of some new heads and wings in Leviathan is an essential advantage in this procedure. The less difference between the two kinds of entropy has been considered as evidence for acceptable modification. 
 
This approach reveals the complexity of the transition between two periodic orbits that exist in the initial and end of the parameter ranges in FHR model. Here, a simple periodic orbit undergoes chaotic changes to arrive at a two-layer one. When the chaotic regime begins, there is a loop in the attractor that can be associated with a wing in the Leviathan. The Leviathan has been discretized by a three vertices graph that two of them are associated with top and down regions in Leviathan and one for the appeared loop.
 When the Markov chain associated with this graph is constructed, its entropy rate can be computed and the comparison between the topological entropy of the dynamical system and the entropy rate of the graph has been done. While the parameters changes whenever there is a noticeable difference between entropy rate and topological entropy, the structure of the Markov chain must be revised. Merging or dividing vertices of the associated graph to the Markov chain are tools to demonstrate the beheading or appearance new heads and wings in the Leviathan.
 
 In the computation of the entropy rate, it is assumed that the walking on the graph of the Markov chain is random. But as it was indicated in the Poincare map of the system there are some regions in which the orbit is periodic and the associated walking is not random at all. 
 Lempel-Ziv of the sequence that is extracted from the simplification of the system is computed for confidence in the randomness of walking. In this paper, the binary version of Lemple-Ziv was applied by simplification of times series into a binary sequence. This quantity certifies us that entropy comparison in all values of parameters is reasonable.
%\nocite{*}
\bibliography{Modified}
%\begin{thebibliography}{99}
%%	\bibitem{10.2307/3649757} 10.2307/3649757
%%	\bibitem{doi:10.1142/S0218127400000840} doi:10.1142/S0218127400000840
%	
%	
%
%\end{thebibliography}
\end{document}